\documentclass[1p]{elsarticle}
\makeatletter
\def\ps@pprintTitle{%
  \let\@oddhead\@empty
  \let\@evenhead\@empty
  \def\@oddfoot{}%
\let\@evenfoot\@oddfoot}
\makeatother
\usepackage{csquotes,microtype,amsmath,mathtools,amsthm}
\usepackage{tikz}
\usepackage{amsfonts,bbm}
\usepackage{hyperref}
\usetikzlibrary{shapes,shadings,patterns}

\usepackage[capitalize]{cleveref}
\crefname{equation}{}{}

\newcommand{\dop}{\mathrm{d}}

\newcommand{\dt}{\dop t}

\newcommand{\ds}{\dop s}
\newcommand{\dW}{\dop W}
\newcommand{\dX}{\dop X}
\newcommand{\dV}{\dop V}

\newcommand{\expo}[1]{\mathrm{e}^{#1}}

\newcommand{\partitionnum}{Q}
\newcommand{\partitionindex}{q}
\newcommand{\vpartitionindex}{r}
\newcommand{\vpartitionnum}{R}
\newcommand{\diffusionnum}{M}
\newcommand{\diffusionindex}{m}
\newcommand{\equaltreeindex}{K}
\newcommand{\numberrepeatedtrees}{v}
\newcommand{\Qindex}{V}
\newcommand{\stagenum}{s}
\newcommand{\mltl}{*}
\newcommand{\treeindex}{k}

\newcommand{\mlt}{\circ}

\newlength{\ml}
  \newcommand{\maketreepartitioned}[1]{
    \tikz [
      grow=up, 				
      level distance=1.3\ml,
      sibling distance=2.7\ml,	
      every node/.style={circle, inner sep=0.05\ml},
      baseline=(current bounding box.base)
    ]
	#1
    ;
  }

  \newcommand{\n}[2]{
	  node[circle, ball color={#1}, inner sep=0.18\ml, label={[yshift=-0\ml]right:{#2}},label={[yshift=+0\ml]left:{\phantom{#2}}}] {}
  }
  \newcommand{\rn}[2]{
	  \node[circle, ball color={#1}, inner sep=0.18\ml, label={[yshift=-0\ml]right:{#2}},label={[yshift=+0\ml]left:{{\phantom{#2}}}}] {}
  }

\newcommand{\treeorig}{
    \maketreepartitioned{\rn{black}{$_{12}$}
    child{\n{black}{$_{21}$}
    }
    child{\n{white}{$_{111}$}
    child{\n{black}{$_{11}$}}
    child{\n{black}{$_{12}$}
    child{\n{black}{$_{21}$}}
    child{\n{black}{$_{21}$}}}}
}}

\newcommand{\treera}{
    \maketreepartitioned{\rn{black}{$_{12}$}
    child[grow=45level distance=1.85\ml]{\n{black}{$_{21}$}}
}}

\newcommand{\treerasemlexp}{
    \maketreepartitioned{\rn{red}{$_{A}$}
    child[grow=45level distance=1.85\ml]{\n{black}{$_{t}$}}
}}

\newcommand{\treena}{
    \maketreepartitioned{\rn{white}{$_{111}$}
    child{\n{black}{$_{11}$}
    }
    child{\n{black}{$_{12}$}
    child{\n{black}{$_{21}$}}
    child{\n{black}{$_{21}$}}}
}}

\newcommand{\treenasemlexp}{
    \maketreepartitioned{\rn{white}{$_{1}$}
    child{\n{black}{$_{0}$}
    }
    child{\n{red}{$_{A}$}
    child{\n{black}{$_{t}$}}
    child{\n{black}{$_{t}$}}}
}}

\newcommand{\treerb}{
    \maketreepartitioned{\rn{black}{$_{12}$}
    child{\n{white}{$_{111}$}
    child{\n{black}{$_{11}$}}
    child{\n{black}{$_{12}$}
    child{\n{black}{$_{21}$}}
    child{\n{black}{$_{21}$}}}}
}}

\newcommand{\treenb}{
    \maketreepartitioned{\rn{black}{$_{21}$}
}}

\newcommand{\treenba}{
    \maketreepartitioned{\rn{black}{$_{12}$}
}}

\newcommand{\treenbb}{
    \maketreepartitioned{\rn{white}{$_{111}$}
}}

\newcommand{\treerc}{
    \maketreepartitioned{\rn{black}{$_{12}$}
    child{\n{black}{$_{21}$}
    }
    child{\n{white}{$_{111}$}
    child{\n{black}{$_{12}$}
    child{\n{black}{$_{21}$}}
    child{\n{black}{$_{21}$}}}}
}}

\newcommand{\treenc}{
    \maketreepartitioned{\rn{black}{$_{11}$}
}}

\newcommand{\treerd}{
    \maketreepartitioned{\rn{black}{$_{12}$}
    child{\n{black}{$_{21}$}}
    child{\n{white}{$_{111}$}
    child{\n{black}{$_{11}$}}
    child{\n{black}{$_{12}$}
    child{\n{black}{$_{21}$}}}}
}}

\newcommand{\treere}{
    \maketreepartitioned{\rn{black}{$_{12}$}
    child{\n{black}{$_{21}$}}
    child{\n{white}{$_{111}$}
    child[grow=45,level distance=1.85\ml]{\n{black}{$_{11}$}}}
}}

\newcommand{\treene}{
    \maketreepartitioned{\rn{black}{$_{12}$}
    child{\n{black}{$_{21}$}}
    child{\n{black}{$_{21}$}}
}}

\newcommand{\texa}{
    \maketreepartitioned{\rn{red}{$_A$}
    child{\n{black}{$_{t}$}
    }
    child{\n{white}{$_{1}$}
    child{\n{black}{$_{0}$}}
    child{\n{red}{$_A$}
      child{\n{black}{$_t$}}
    child{\n{black}{$_t$}}}}
}}

\newtheorem{example}{Example}
\newtheorem{definition}{Definition}
\newtheorem{theorem}{Theorem}
\newtheorem{lemma}{Lemma}

\newtheorem{corollary}{Corollary}

\colorlet{added}{purple}
\colorlet{changed}{teal}

\begin{document}
\title{B-series for SDEs with application to exponential integrators for non-autonomous semi-linear problems}
\author[Hawassa]{Alemayehu Adugna Arara}
\ead{alemayehuadugna@hu.edu.et}
\author[SDU]{Kristian Debrabant}
\ead{debrabant@imada.sdu.dk}
\author[NTNU]{Anne Kv\ae{}rn\o{}}
\ead{anne.kvarno@ntnu.no}
\address[Hawassa]{Hawassa University, Department of Mathematics, P.O. Box 05, Hawassa, Ethiopia}
\address[SDU]{University of Southern Denmark, Department of Mathematics and Computer Science, 5230 Odense M, Denmark}
\address[NTNU]{Norwegian University of Science and Technology - NTNU, Department of Mathematical Sciences, NO-7491 Trondheim, Norway}
\begin{abstract}
  In this paper a set of previous general results for the development of B--series for a broad class of stochastic differential equations has been collected. The applicability of these results is demonstrated by the derivation of B--series for non-autonomous semi-linear SDEs and exponential Runge-Kutta methods applied to this class of SDEs, which is a significant generalization of existing theory on such methods.
\end{abstract}
\maketitle
\begin{center}
\textit{Dedicated to Brynjulf Owren and Hans Munthe-Kaas on the occasion of their 60th birthdays}
\end{center}
\section{Introduction}

The idea of B--series, as described here, is simply about describing the solution of a differential equation (or approximating numerical scheme) as
a formal series in time (or in the step size){, allowing then to represent the local error, i.e.\ the difference between the exact solution and the numerical approximation after one time step, as a B-series as well and thereby easily analyse the order of convergence of the numerical scheme}. These series were initially derived for ordinary differential equations (ODEs) by Butcher \cite{butcher63cft}, and the idea was further extended by Hairer and Wanner in \cite{hairer73mmm}, see \cite{hairer10gni} and references therein. Owren and Munthe-Kaas have both given several contributions to the development of B--series, in particular in the context of geometrical integration, see  \cite{ munthekaas95lbt,  berland05bsa, owren06ocf, celledoni10epi, munthekaas13opl, munthe-kaas16abs}. B--series for stochastic differential equations (SDEs) have been developed in different contexts by several authors, \cite{burrage96hso, burrage99thesis,
  burrage00oco, komori97rta, komori07mcr, roessler04ste, roessler06rta}.  A rather general framework
for developing B--series for SDEs was developed and extended in \cite{debrabant08bsa, debrabant11cos, anmarkrud18goc}.
This approach is independent on whether the SDE is It\^{o} or Stratonovich and whether weak or strong convergence is considered.

The aim of this paper is to collect the main tools from this approach, which is done in section \cref{sec:bseries}, and then apply the theory to non-autonomous semi-linear SDEs and a class of exponential Runge--Kutta (ERK) methods. A similar work has very recently been published by Yang et.al.\ \cite{yang22anc}, however, their results only apply to Stratonovich SDEs with a commutative linear part. The results which will be presented in the current paper do not have these limitations.

The outline is as follows: The main ideas and results on stochastic B--series for autonomous SDEs are presented in \cref{sec:bseries}. In \cref{sec:non-autonomous} these are used for a generalisation to non-autonomous SDEs, and non-autonomous semi-linear SDEs are discussed in \cref{sec:semilinearSDEt}.

\section{Stochastic B--series: Definitions and main results} \label{sec:bseries}
In this section the basic definitions and results to be used in the later sections are
summarised.

Consider a system of stochastic partitioned differential equations with
$\partitionnum$ partitions and $\diffusionnum$ diffusion terms that potentially are split,
\begin{equation}
X^{(\partitionindex)}(t)=x_0^{(\partitionindex)}+\sum_{\diffusionindex=0}^\diffusionnum\int_0^t
\big(\sum_{\vpartitionindex=1}^{\vpartitionnum_{\diffusionindex,\partitionindex}}
g^{(\partitionindex,\vpartitionindex)}_{\diffusionindex}(X^{(1)}(s),X^{(2)}(s),\dots,X^{(\partitionnum)}(s))\big)\star
\dW_\diffusionindex(s)\label{eq:SDEIntegralform}
\end{equation}
for $\partitionindex=1,\ldots,\partitionnum$, for which we will also use the abbreviated form
\begin{equation}
\dop X^{(\partitionindex)}(t)=\sum_{\diffusionindex=0}^\diffusionnum
\big(\sum_{\vpartitionindex=1}^{\vpartitionnum_{\diffusionindex,\partitionindex}}
g^{(\partitionindex,\vpartitionindex)}_{\diffusionindex}(X^{(1)}(t),X^{(2)}(t),\dots,X^{(\partitionnum)}(t))\big)\star
\dW_\diffusionindex(t),
\quad X^{(q)}(0)=x^{(q)}_0. \label{eq:SDE}
\end{equation}
To simplify the notation the deterministic terms are represented by $\diffusionindex=0$, such that $\dW_0(s)=\ds$, while $W_\diffusionindex$, $\diffusionindex=1,\dots,\diffusionnum$, denote one-dimensional and pairwise independent Wiener processes.
The integrals w.r.t.\ the Wiener processes are interpreted as either {It\^o} integrals,
$\star \dW_\diffusionindex(s)=\dW_\diffusionindex(s)$, or Stratonovich integrals, $\star \dW_\diffusionindex(s)=\circ \dW_\diffusionindex(s)$.
We also define the vector of initial values, $x_0 = [{x^{(1)}_{0}}^\top, {x^{(2)}_{0}}^\top, \dots, {x^{(\partitionnum)}_{0}}^\top]^\top$.
Furthermore, we assume that the coefficients $g^{(\partitionindex,\vpartitionindex)}_{\diffusionindex}: \mathbb{R}^{d_1}\times \ldots \times \mathbb{R}^{d_\partitionnum} \to \mathbb{R}^{d_q}$ are sufficiently smooth,
and that the conditions of the existence and uniqueness theorem \cite{oksendal03sde} are satisfied.
The systems are considered to be autonomous.
If $\vpartitionnum_{\diffusionindex,\partitionindex}=1$ for $\diffusionindex=0,\dots,\diffusionnum$ and $\partitionindex=1,\dots,\partitionnum$ then the splitting is  horizontal, if $\partitionnum=1$ the splitting is vertical, see
\cite{ascher97ier}. In those cases, the corresponding indices may be omitted.

{
\begin{example} \label{ex:example}
Consider the equation system
    \begin{align*}
      \dot{r} &= v, \\
      \dot{v} &= f_d(r,t) - \alpha(t)v + f_s(r,t)\beta(t)
    \end{align*}
  which describes the evolution of a particle with unit mass, coordinate $r(t)$ and velocity
  $v(t)$. The particle is affected by three forces: an external force $f_d(t,r)$, a time-dependent
  friction force $\alpha(t)\geq 0$ and a thermal white noise term with a time and space dependent intensity. This can be rewritten as a proper SDE of the form \cref{eq:SDE},
  \begin{subequations} \label{eq:langevin_vertical}
    \begin{align}
      \dop\underbrace{\begin{pmatrix} R(t) \\ V(t) \end{pmatrix}}_{X^{(1)}(t)} &=
          \underbrace{\begin{pmatrix} 0 \\ f_d(R(t),t) \end{pmatrix}}_{g_0^{(1,1)}} \dt +
          \underbrace{\begin{pmatrix} V(t) \\ -\alpha(t) V(t) \end{pmatrix}}_{g_0^{(1,2)}} \dt +
          \underbrace{\begin{pmatrix} 0 \\ f_s(R(t),V(t),t) \end{pmatrix}}_{g_1^{(1,1)}} \dW_1(t) \\
      \dop\underbrace{t}_{X^{(2)}(t)}&= \underbrace{1}_{g_0^{(2,1)}}\dt.
    \end{align}
  \end{subequations}
  This is of course only one of
  many possible splittings of this problem, but it is the one which will be used for demonstration purposes throughout this
  paper. Other splittings have been discussed in  \cite{anmarkrud18goc}, Example
  1.1.
\end{example}
}

Partitioned system{s} appea{r} natural{ly} in many situations, e.g.\ for Hamiltonian problems. {Note that \cref{eq:SDE} also covers the case of matrix-valued SDEs, in which case $X^{(q)}$ would correspond to column $q$ of the solution.}

In this
paper, horizontal splitting will be used to derive B--series for non-autonomous problems, vertical
splitting will be used for semi-linear problems.

The main idea is to express the exact solutions $X^{(\partitionindex)}(h)$  as B--series $B^{(\partitionindex)}(\phi,x_0;h)$:
\begin{align*}
B^{(\partitionindex)}(\phi,x_0;h)=\sum_{\tau \in T_\partitionindex}\alpha(\tau)\cdot \phi(\tau)(h)\cdot F(\tau)(x_0),
\end{align*}
where $T_\partitionindex$ is the set of shaped, colored, rooted trees as defined below.
The terms $\alpha(\tau)$ are combinatoric terms. The elementary weight functions $\phi(\tau)(h)$ are stochastic integrals or random variables, and $F(\tau)(x_0)$ are the elementary differentials. To simplify the presentation, we assume that all elementary differentials exist and all considered B--series converge. Otherwise, one has to consider truncated B-series and discuss the remainder term \cite{roessler04ste}.

The following theory on partitioned SDEs is taken from \cite{anmarkrud18goc}, where also all proofs can be found.

\begin{definition}[Trees and combinatorial coefficients]\label{def:trees}
The set of shaped, rooted trees
\begin{align*}
  T =T_1 \cup T_2 \cup \dots \cup T_\partitionnum
\end{align*}
where
\begin{align*}
  T_{\partitionindex,\diffusionindex} &= T_{\partitionindex,0,\diffusionindex} \cup
                                        T_{\partitionindex,1,\diffusionindex} \cup \dots \cup T_{\partitionindex,
                                                         \vpartitionnum_{\diffusionindex, \partitionindex},\diffusionindex}
  & \text{and} \qquad
  T_\partitionindex &= \{\emptyset_\partitionindex\} \cup T_{\partitionindex,0} \cup T_{\partitionindex,1} \cup \dots \cup T_{\partitionindex,\diffusionnum}
\end{align*}
for $\partitionindex=1,\ldots,\partitionnum$ and $\diffusionindex=0,\dots,\diffusionnum$ is recursively defined as follows for $\vpartitionindex=1,\dots,\vpartitionnum_{\diffusionindex,\partitionindex}$:
\begin{enumerate}
  \item The graph
    $\bullet_{\partitionindex,\vpartitionindex,\diffusionindex}$ with only one vertex of shape
    $(\partitionindex,\vpartitionindex)$ and color $\diffusionindex$ belongs to $T_{\partitionindex,\vpartitionindex,\diffusionindex}$.
  \item If $\tau_1,\tau_2,\dots,\tau_{\kappa} \in
    T\setminus\{\emptyset_1,\dots,\emptyset_{\partitionnum}\}$, then also
    $[\tau_1,\tau_2,\dots,\tau_{\kappa}]_{\partitionindex,\vpartitionindex,\diffusionindex}\in T_{\partitionindex,\vpartitionindex, \diffusionindex}$, where
    $[\tau_1,\tau_2,\dots,\tau_{\kappa}]_{\partitionindex,\vpartitionindex,\diffusionindex}$
    denotes the tree formed by joining the subtrees
    $\tau_1,\tau_2,\dots,\tau_{\kappa}$ each by a single branch to a common root of
    shape $(\partitionindex,\vpartitionindex)$ and color $\diffusionindex$.
\end{enumerate}
Further, we define $\alpha(\tau)$ as
\begin{align*}
  \alpha(\emptyset_\partitionindex) = 1, \;
  \alpha(\bullet_{\partitionindex,\vpartitionindex,\diffusionindex}) = 1, \;
  \alpha(\tau=[\tau_1, \dots, \tau_{\kappa}]_{\partitionindex,\vpartitionindex,\diffusionindex}) = \frac{1}{\numberrepeatedtrees_1!\numberrepeatedtrees_2!\dots \numberrepeatedtrees_\equaltreeindex!}\prod_{\treeindex=1}^{\kappa}\alpha(\tau_\treeindex),
\end{align*}
where $\numberrepeatedtrees_1,\numberrepeatedtrees_2,\dots,\numberrepeatedtrees_\equaltreeindex$ count equal trees among $\tau_1,\tau_2,\dots,\tau_{\kappa}$.
\end{definition}

In addition to the bracket notation used in \cref{def:trees}, rooted trees can be illustrated as graphs, see \cref{ex:thmexact}, where vertices $\bullet_{\partitionindex,\vpartitionindex,0}$ with color $\diffusionindex=0$ are  represented as black nodes (so-called deterministic nodes) with indices $\partitionindex,\vpartitionindex$, while vertices $\bullet_{\partitionindex,\vpartitionindex,\diffusionindex}$ with color $\diffusionindex>0$ are represented by white nodes (so-called stochastic nodes) with indices $\partitionindex,\vpartitionindex,\diffusionindex$.

\begin{definition}[Elementary differentials] For a tree $\tau \in T$ the elementary differential is a mapping $F(\tau)$: $\mathbb{R}^{d_1}\times \ldots \times \mathbb{R}^{d_\partitionnum} \to \mathbb{R}^{d}$ defined recursively by
\begin{enumerate}
  \item $F(\emptyset_\partitionindex)(x_0) = x_{0}^{(\partitionindex)}$, %$\emptyset_\partitionindex \in T_\partitionindex$,
  \item
    $F(\bullet_{\partitionindex,\vpartitionindex,\diffusionindex})(x_0)=g^{(\partitionindex,\vpartitionindex)}_{\diffusionindex}(x_0)$,
  \item {i}f $\tau =
    [\tau_1,\tau_2,\dots,\tau_{\kappa}]_{\partitionindex,\vpartitionindex,\diffusionindex}\in T_{\partitionindex,\vpartitionindex,\diffusionindex}$, then
    \begin{equation*}
      F(\tau)(x_0)=(D_{\tau}g^{(\partitionindex,\vpartitionindex)}_{\diffusionindex})(x_0)(F(\tau_1)(x_0),F(\tau_2)(x_0),\dots,F(\tau_{\kappa})(x_0))
    \end{equation*}
    where $D_{\tau}=\frac{\partial^{\kappa}}{\partial x^{(\partitionindex_1)}\dots\partial x^{(\partitionindex_\kappa)}}$
    for $(\partitionindex_\treeindex,\vpartitionindex_\treeindex)$ being the shape of $\tau_\treeindex$, $\treeindex=1,\dots,\kappa$.
\end{enumerate}\label{def:elementarydifferentials}
\end{definition}
Fundamental for this work is the following lemma which says that if {some functions} ${V}^{(\partitionindex)}(h)$ can be written as a B-series, then $f({V}^{(1)}(h),\dots,{V}^{(\partitionnum)}(h))$ can also be written as a B-series.
This is a trivial extension of the lemma found in \cite{debrabant08bsa}.
\begin{lemma}\label{lem:function_of_b-series}
If ${V}^{(\partitionindex)}(h) = B^{(\partitionindex)}(\phi,x_0;h)$, $\partitionindex=1,\dots,\partitionnum$, are some B-series with $\phi(\emptyset_q)\equiv 1$ and $f \in C^{\infty}(\mathbb{R}^{d_1}\times \ldots \times \mathbb{R}^{d_\partitionnum}, \mathbb{R}^{d})$, then $f({V}^{(1)}(h),\dots,{V}^{(\partitionnum)}(h))$ can be written as a formal series of the form
\begin{equation}
  f({V}^{(1)}(h),\dots,{V}^{(\partitionnum)}(h)) = \sum_{u\in U_f} \beta(u)\cdot \psi_{\phi}(u)(h)\cdot G(u)(x_0)
\end{equation}
where
\begin{enumerate}
  \item $U_f$ is a set of trees derived from T, by
  \begin{enumerate}
    \item $\bullet_f \in U_f$,
    \item if $\tau_1,\tau_2,\dots,\tau_{\kappa} \in T\setminus\{\emptyset_1,\dots,\emptyset_{\partitionnum}\}$, then $[\tau_1,\tau_2,\dots,\tau_{\kappa}]_f \in U_f$.
  \end{enumerate}
  \item $G(\bullet_f)(x_0)=f(x_0)$ and
    \[G(\tau=[\tau_1,\tau_2,\dots,\tau_{\kappa}]_f)(x_0)=({D_{\tau}} f)(x_0)(F(\tau_1)(x_0),\dots,F(\tau_{\kappa})(x_0)),\]
    where $D_\tau$ is defined in analogy to \cref{def:elementarydifferentials}.
  \item $\beta(\bullet_f)=1$ and \[\beta([\tau_1,\tau_2,\dots,\tau_{\kappa}]_f)=\frac{1}{\numberrepeatedtrees_1!\numberrepeatedtrees_2!\cdots \numberrepeatedtrees_\equaltreeindex!}\prod_{\treeindex=1}^{\kappa}\alpha(\tau_\treeindex),\] where $\numberrepeatedtrees_1,\numberrepeatedtrees_2,\dots,\numberrepeatedtrees_\equaltreeindex$ count equal trees among $\tau_1,\tau_2,\dots,\tau_{\kappa}$.
  \item $\psi_{\phi}(\bullet_f)(h){=} 1$ and  $\psi_{\phi}([\tau_1,\tau_2,\dots,\tau_{\kappa}]_f)(h)=\prod_{\treeindex=1}^{\kappa}\phi{(\tau_\treeindex)}(h)$.
\end{enumerate}
\end{lemma}
By the use of this lemma, the following resul{t was} proved in \cite{anmarkrud18goc}.
\begin{theorem}\label{th:b-seriesexact}The exact solutions $X^{(\partitionindex)}({h})$ of \cref{eq:SDEIntegralform}, $\partitionindex=1,\dots,\partitionnum$, can be written as B-series $B^{(\partitionindex)}(\varphi,x_0;h)$ with
\begin{gather*}
  \varphi(\emptyset_\partitionindex)(h){=} 1, \quad \varphi(\bullet_{\partitionindex,\vpartitionindex,\diffusionindex})(h)={\int_0^h\dW_\diffusionindex(s)=}W_\diffusionindex(h), \\
  \varphi([\tau_1,\tau_2,\dots,\tau_\kappa]_{\partitionindex,\vpartitionindex,\diffusionindex})(h)=\int_0^h\prod_{j=1}^\kappa\varphi(\tau_j(s)) \star \dW_\diffusionindex(s), \\
  \text{for all }[\tau_1,\tau_2,\dots,\tau_\kappa]_{\partitionindex,\vpartitionindex,\diffusionindex} \in T_{\partitionindex,\diffusionindex}, \;
  \vpartitionindex=1,\dots,\vpartitionnum_{\diffusionindex,\partitionindex},\;
  \partitionindex=1,\dots,\partitionnum, \; \diffusionindex=0,1,\dots,\diffusionnum.
\end{gather*}
\label{thm:exact}
\end{theorem}
{\Cref{th:b-seriesexact} represents the solution of \cref{eq:SDEIntegralform} as B-series. By also deriving a B-series representation of the numerical approximation, one can represent the local error, i.e.\ their difference, as well as B-series. Comparing the weight functions in the B-series representation of the exact solution and numerical approximation allows thus to analyze the (strong or mean-square) consistency of a numerical method. For weak approximations, where function evaluations of the exact solution and numerical approximation need to be compared, additionally \cref{lem:function_of_b-series} is used.}

{To characterize the contribution of (the weight function of) a tree to the order of the error, the following definition is needed.
\begin{definition}[Tree order]
The order $\rho$ of a tree $\tau\in T$ is defined by
\begin{gather*}
  \rho(\emptyset_\partitionindex)= 1, \quad\rho(\bullet_{\partitionindex,\vpartitionindex,\diffusionindex})=\begin{cases}
    1&\text{for }\diffusionindex=0,\\
    \frac12&otherwise,
  \end{cases}\\
  \rho([\tau_1,\tau_2,\dots,\tau_\kappa]_{\partitionindex,\vpartitionindex,\diffusionindex})=\sum_{k=1}^\kappa\rho(\tau_k)+\begin{cases}
    1&\text{for }\diffusionindex=0,\\
    \frac12&otherwise.
  \end{cases}
\end{gather*}
\end{definition}}
\begin{example} \label{ex:thmexact}
  {Let}
  \[
    \tau = \treeorig = [[[\bullet_{2,1,0}^2]_{1,2,0},\bullet_{1,1,0}]_{1,1,1},\bullet_{2,1,0}]_{1,2,0}{,}
  \]
  {in which vertices $\bullet_{\partitionindex,\vpartitionindex,0}$ with color $\diffusionindex=0$ are represented as black nodes with indices $\partitionindex,\vpartitionindex$, while vertices $\bullet_{\partitionindex,\vpartitionindex,1}$ with color $\diffusionindex=1$ are represented by white nodes with indices $\partitionindex,\vpartitionindex,1$.}
  {T}he corresponding terms are
  \begin{align*}
    \alpha(\tau) &= 1/2,  \\
    {\rho(\tau)}&{=6\frac12,}\\
    F(\tau) &= D_{1,2}g_{{0}}^{({1,2})}\Bigg(D_{1,1}g_1^{(1,1)}\Big(D_{2,2}g_{{0}}^{({1,2})}\big(g_{{0}}^{({2,1})},g_{{0}}^{({2,1})}\big),g_0^{(1,{1})}\Big),g_{{0}}^{({2,1})}\Bigg), \\
    \varphi(\tau) &= \int_0^h  \int_0^{s_1}  \int_0^{s_2} s_3^2\ds_3 s_2 \star \dW_1(s_2)  s_1 \ds_1{=\frac13\int_0^h  \int_0^{s_1}  s_2^4 \star \dW_1(s_2)  s_1 \ds_1}
  \end{align*}
  where $D_{ij}=\frac{\partial^2}{\partial x^{(i)}\partial x^{(j)}}$.
  {How the elementary differentials would look like for the SDE introduced in \cref{ex:example} is presented in \ref{sec:appA} (the weight function $\varphi(\tau)$ and coefficient $\alpha$ only depend on the tree, not on the concrete SDE).}
\end{example}
Stochastic B--series evaluated in a B--series can (similar to the deterministic case) be written as a B--series. This result is useful for the construction of implicit Taylor methods {(see e.\,g.\ \cite{kloeden99nso,tian01itm,milstein02nmf,debrabant11bsa,debrabant11cos})}, splitting methods and composition methods. {Note that for implicit methods, it might be necessary to use truncated random variables to ensure the solvability of the implicit equations, see \cite{milstein02nmf}.}

Let $ST(\tau)$ be the set of all possible subtrees of $\tau$ together with the corresponding  remainder multiset $\omega$, that is, for each $\tau\in T\setminus\{\emptyset_1,\dots,\emptyset_{\partitionnum}\}$ we have
\begin{align*}
ST(\bullet_{\partitionindex,\vpartitionindex,\diffusionindex}) &= \{ (\emptyset_\partitionindex,\{\bullet_{\partitionindex,\vpartitionindex,\diffusionindex}\}), (\bullet_{\partitionindex,\vpartitionindex,\diffusionindex},\{\emptyset_\partitionindex\}) \}, \\
ST(\tau =[\tau_1,\dotsc,\tau_{\kappa}]_{\partitionindex,\vpartitionindex,\diffusionindex}) &=
\{(\vartheta,\omega)\;:\; \vartheta=[\vartheta_1,\dotsc,\vartheta_{\kappa}]_{\partitionindex,\vpartitionindex,\diffusionindex}, \quad
   \omega=\{\omega_1,\dotsc,\omega_{\kappa}\}, \\ &\qquad
   (\vartheta_i,\omega_i) \in ST(\tau_i), \quad i=1,\dots,\kappa \} \cup
   (\emptyset,\{\tau\}).
\end{align*}
We also have to take care of possible equal terms in the formula presented below. This is done as follows:
For a given triple $(\tau,\vartheta,\omega)$ write first $\vartheta=[\vartheta_1,\dotsc,\vartheta_{\kappa_{\vartheta}}]_{\partitionindex,\vpartitionindex,\diffusionindex} = [\bar{\vartheta}_1^{{\numberrepeatedtrees}_1},\dotsc,\bar{\vartheta}_{\equaltreeindex}^{{\numberrepeatedtrees}_\equaltreeindex}]_{\partitionindex,\vpartitionindex,\diffusionindex}$, where the latter only expresses that $\vartheta$ is composed by $\equaltreeindex$ different nonempty trees, each appearing $\numberrepeatedtrees_i$ times, hence $\sum_{i=1}^\equaltreeindex \numberrepeatedtrees_i=\kappa_{\vartheta}$.
Let $\tau=[\tau_1,\dotsc,\tau_{\kappa}]_{\partitionindex,\vpartitionindex,\diffusionindex}$. For $i=1,\dots,\equaltreeindex$, each $\bar{\vartheta}_i$ is a subtree of some of the $\tau_j$'s, with corresponding remainder multisets $\omega_j$. Assume that there are exactly $p_i$ different such triples $(\bar{\tau}_{ik},\bar{\vartheta}_i,\bar{\omega}_{ik})$ each appearing exactly $\numberrepeatedtrees_{ik}$ times so that $\sum_{k=1}^{p_i}\numberrepeatedtrees_{ik}=\numberrepeatedtrees_i$.
Finally, let $\bar\delta_k\in\omega$ be the distinct trees with multiplicity $s_k$, $k=1,\dots,p_0$, of the remainder multiset which are directly connected to the root of $\tau$.
Then, $\tau$ can be written as
\begin{equation}
\label{bsa2:eq:alltrees}
\tau=[\bar\delta_1^{s_1},\dotsc, \bar\delta_{p_0}^{s_{p_0}},
    \bar{\tau}_{11}^{\numberrepeatedtrees_{11}},\dotsc,\bar{\tau}_{1p_1}^{\numberrepeatedtrees_{1p_1}} ,\dotsc,
    \bar{\tau}_{\equaltreeindex1}^{\numberrepeatedtrees_{\equaltreeindex1}},\dotsc,\bar{\tau}_{{\equaltreeindex}p_{\equaltreeindex}}^{\numberrepeatedtrees_{{\equaltreeindex}p_{\equaltreeindex}}}]_{\partitionindex,\vpartitionindex,\diffusionindex}=
    [ \bar{\tau}_1^{R_1},\dotsc,\bar{\tau}_{\Qindex}^{R_\Qindex}]_{\partitionindex,\vpartitionindex,\diffusionindex},
\end{equation}
where the rightmost expression above indicates that $\tau$ is composed by $\Qindex$ different trees $\bar{\tau}_1,\dots,\bar{\tau}_\Qindex$, with $\bar{\tau}_i$ appearing $R_i$ times, $i=1,\dots,\Qindex$.

With these definitions, we can state the following theorem, {whose} proof is similar to the one for the nonpartitioned case given in  \cite{debrabant11cos}, see also \cite{debrabant11bsa}:
\begin{theorem}[Composition of B--series] \label{thm:comp}
Let $\phi_x(\emptyset_\partitionindex)\equiv 1$ for $\partitionindex=1,\dots,\partitionnum$. Then for $\partitionindex=1,\dots,\partitionnum$, the B--series $B^{(\partitionindex)}(\phi_y,\cdot;h)$ evaluated at
\[B(\phi_x,x_0;h):=[B^{(1)}(\phi_x,x_0;h)^\top,\dots,B^{(\partitionnum)}(\phi_x,x_0;h)^\top]^\top\] is again a B--series,
\[ B^{(\partitionindex)}(\phi_y,B(\phi_x,x_0;h);h)=B^{(\partitionindex)}(\phi_x\mlt \phi_y,x_0;h),\]
where
\[ \left(\phi_x\mlt \phi_y \right)(\tau) = \sum_{(\vartheta,\omega)\in ST(\tau)} \gamma(\tau,\vartheta,\omega)\left( \phi_y(\vartheta)\prod_{\delta \in \omega} \phi_x(\delta)\right)
\]
with $\gamma(\emptyset_\partitionindex,\emptyset_\partitionindex,\{\emptyset_\partitionindex\})=1$ and
\begin{equation*}%\label{eq:gamma}
\gamma(\tau,\vartheta,\omega)=\frac{R_1!\dotsm R_\Qindex!}{s_1!\dotsm s_{p_0}! \numberrepeatedtrees_{11}!\dotsm \numberrepeatedtrees_{{\equaltreeindex}p_{\equaltreeindex}}!} \prod_{i=1}^{\kappa_{\vartheta}} \gamma(\tau_i,\vartheta_i,\omega_i).
\end{equation*}
for $\tau\neq\emptyset_\partitionindex$.
\end{theorem}

Finally, the following lemma, proved in \cite{debrabant11bsa}, is the key to derive  B--series for
exponential integrators for SDEs in \cref{sec:semilinearSDEt}.
\begin{lemma} \label{bsa2:lem:Bfprime}
If $\phi_x(\emptyset_\partitionindex)\equiv0$ for $\partitionindex=1,\dots,\partitionnum$, then we have for $\partitionindex=1,\dots,\partitionnum$
\[
\partial_2B^{(\partitionindex)}(\phi_{y},{x_{0};h})B(\phi_{x},{x_{0};h})
=B^{(\partitionindex)}(\phi_{x}\mltl\phi_{y},{x_{0};h}),
\]
where the bi-linear operator $\mltl$ is given by
\begin{equation} \label{bsa2:eq:mltl}
(\phi_x\mltl\phi_y)(\tau)=\begin{cases}
  0
  &if \; \tau=\emptyset,\\
  \sum\limits_{(\vartheta,\{\delta\})\in SP(\tau)}\gamma(\tau,\vartheta,\{\delta\})\cdot
  \phi_{y}(\vartheta)\phi_{x}
  (\delta) &\text{otherwise},
  \end{cases}
\end{equation}
with
\[
SP(\tau)=\{(\vartheta,\omega)\in ST(\tau):~\omega\text{ contains
  exactly one element } \delta\}.
\]
\end{lemma}

\begin{example}\label{ex:sp}
  Consider the tree $\tau$ from \cref{ex:thmexact},
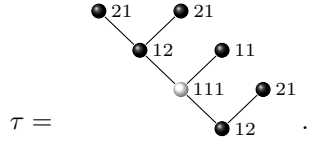
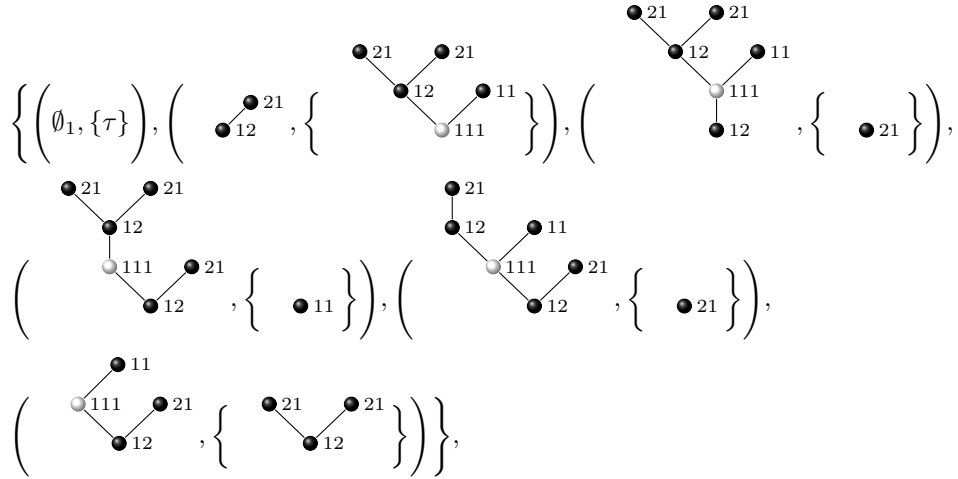
  \[ \tau = \treeorig. \]
  Then
  \begin{align*}
    & SP(\tau)  = \\
    & \Bigg\{ \Bigg( \emptyset_1, \{\tau\} \Bigg),
      \Bigg(\treera, \bigg\{\treena\bigg\}\Bigg),
      \Bigg( \treerb, \bigg\{\treenb \bigg\} \Bigg),\\
      &\Bigg(\treerc, \bigg\{\treenc\bigg\} \Bigg),
      \Bigg(\treerd, \bigg\{\treenb \bigg\} \Bigg), \\
      &\Bigg(\treere, \bigg\{\treene\bigg\}\Bigg)
       \Bigg\}{,}
    \end{align*}{ which will be further used in \cref{ex:exptree}.}
\end{example}

For general SDEs, the number of terms in the B--series is rather overwhelming, so whenever applied to a particular class of SDEs, the main issue is to identify all trees for which the elementary differentials automatically become zero, and remove those from the definition of trees for this class. This idea was demonstrated for several problems in \cite{anmarkrud18goc} and will again be applied here.

\section{Non-autonomous SDEs}\label{sec:non-autonomous}
B-series are usually expressed for autonomous versions of the equations, with the assumption, which
also will be used here, that a non-autonomous system easily can be converted to an autonomous one.
This is of particular interest in the rather common case of SDEs with
additive but time-dependent noise, for example:

\begin{equation}
  \label{eq:}
  \dX(t) = g_0(X(t){)}\dt + \sum_{\diffusionindex = 1}^\diffusionnum g_\diffusionindex(t) \star \dW_\diffusionindex(t).
\end{equation}

In this section, we will derive B--series for SDEs for which also some or all of the Wiener processes are allowed
to enter explicitly into the coefficient functions,
\begin{equation} \label{eq:SDE_nonautonomous}
  \dop X(t)= \sum_{m=0}^{M}\sum_{\vpartitionindex=1}^{\vpartitionnum_\diffusionindex}g_{m}^{(\vpartitionindex)}\big(X(t),W(t)\big) \star \dW_m(t),\quad X(0)=x_{0},
\end{equation}
where $W(t)=[W_0(t), W_1(t), \dotsc, W_l(t)]^\top$ with $l\leq m$ and we allow for vertical splitting.

\begin{example}[Lawson - or the integrating factor method]
  Consider the semi-linear equation
  \begin{equation} \label{eq:multilin}
    \dX(t) = \sum_{\diffusionindex=0}^\diffusionnum \big( A_\diffusionindex X(t) + g_\diffusionindex(X(t) \big)\star \dW_\diffusionindex(t)
  \end{equation}
  with constant, commutative linear terms, that is
\[ [A_{\diffusionindex_1},A_{\diffusionindex_2}] = A_{\diffusionindex_1}A_{\diffusionindex_2}-A_{\diffusionindex_2}A_{\diffusionindex_1} = 0, \qquad {\diffusionindex_1},{\diffusionindex_2} = 1,\dotsc,\diffusionnum.\]
Using the Lawson transformation
\[
  V(t) = \expo{-L(t)}X(t), \qquad L(t) = \big(A_0 - \gamma^{\star} \sum_{\diffusionindex=1}^{\diffusionnum}{A_\diffusionindex^2})t + \sum_{\diffusionindex=1}^\diffusionnum
  A_m W_m(t)
\]
with $\gamma^\star=1/2$ in the It\^{o} case and ${\gamma^\star=}0$ in the Stratonovich case,
the transformed system becomes
\[
  \dV(t) = \sum_{\diffusionindex=0}^\diffusionnum \hat{g}_\diffusionindex(V(t),W(t))\star \dW_\diffusionindex(t), \qquad
  \hat{g}_\diffusionindex(x,W(t)) = \expo{-L(t)}g_\diffusionindex(\expo{L(t)}{x})
\]
which is a non-autonomous SDE of the form \cref{eq:SDE_nonautonomous}. See \cite{debrabant21rkl} and references therein.
\end{example}

Following the well known approach from ODEs, let  $\tilde{X}^{(\diffusionindex+2)}(t) =
{W_{\diffusionindex}(t)}$, $\diffusionindex = 0,\dotsc,l$, such that the SDE can be written
as a horizontal (and vertical) split system:
\begin{align*}
  \dop\tilde{X}^{(1)}(t) &= \sum_{\diffusionindex=0}^{\diffusionnum}\sum_{\vpartitionindex=1}^{\vpartitionnum_\diffusionindex}g^{(1,\vpartitionindex)}_{\diffusionindex}(\tilde{X}(t)) \star
  \dW_\diffusionindex(t), & \tilde{X}^{(1)}(0)&=x_{0}, \\
  \dop\tilde{X}^{(2)}(t) &= g_{0}^{(2,1)}(\tilde{X}(t))\star\dW_{0}(t) = 1 \cdot \dW_{0}(t), & \tilde{X}^{(2)}(0) &= t_0, \\
  & \vdots & & \vdots \\
  \dop\tilde{X}^{(l+2)}(t) &= g_{l}^{(l+2,1)}(\tilde{X}(t))\star\dW_{l}(t) = 1  \cdot \dW_{l}(t),	& \tilde{X}^{(l+2)}(0) &= W_l(t_0),
\end{align*}
where $\tilde{X}(t) = [X(t+t_0)^{\top},W_0(t+t_0),\dotsc,W_l(t+t_0)]^\top$ and $g^{(1,\vpartitionindex)}_{\diffusionindex}=g^{(\vpartitionindex)}_{\diffusionindex}$.
To simplify the notation, we will in the following write $\bullet_{\vpartitionindex,\diffusionindex}$ instead of $\bullet_{1,\vpartitionindex,\diffusionindex}$, and
denote the nodes corresponding
to the $W_i$'s by $\bullet_{W_i}$ or also simply $\bullet_t$ in case of $i=0$.
For their elementary differentials it holds that $F(\bullet_{W_i})(X(t))=1$, and all further differentials of these are 0, and their corresponding trees can be omitted from the set of trees.
With these considerations we conclude:

\begin{corollary}\label{cor:autonomous}
  The solution $X(t)$ of \cref{eq:SDE_nonautonomous} can be written as a B--series for which:
  \begin{enumerate}
    \item The set of trees  $T$  is defined by
      \begin{gather*}
	\emptyset \in T \text{ and } \bullet_{\vpartitionindex,\diffusionindex}
	\in T , \\
	\tau = [\tau_1,\dotsc,\tau_{\kappa}]_{\vpartitionindex,\diffusionindex} \in T \text{ for all }
	\tau_1, \dotsc, \tau_{\kappa} \in T\cup\{\bullet_{W_0},\dots,\bullet_{W_\diffusionnum}\},
      \end{gather*}
      where $\vpartitionindex=1,\dots,\vpartitionnum_\diffusionindex$, $\diffusionindex=0,\dotsc,\diffusionnum$.
    \item The elementary differentials are given by
      \begin{gather*}
	F(\bullet_{\vpartitionindex,\diffusionindex})(\tilde{x}_0) = g_\diffusionindex^{(\vpartitionindex)}(\tilde{x}_0), \qquad
	F(\bullet_{W_\diffusionindex})(\tilde{x}_0) = 1, \\
	F(\tau = [\tau_1,\tau_2,\dotsm,\tau_{\kappa}]_{\vpartitionindex,\diffusionindex})(\tilde{x}_0)
	= \left(D_{\tau}g_{\diffusionindex}^{(\vpartitionindex)}(\tilde{x}_0)\right)
	\left( F(\tau_1)(\tilde{x}_0),\dotsc,F(\tau_{\kappa})(\tilde{x}_0) \right).
    \end{gather*}
\end{enumerate}
\end{corollary}

\section{Semi-linear SDEs and exponential integrators } \label{sec:semilinearSDEt}
Semi-linear problems have over the years  been given quite some attention, and several
exponential methods have been developed both for ODEs and SDEs to exploit the special structure of
such systems. B--series for semi-linear ODEs have been discussed e.g.\  in \cite{berland05bsa,
  hochbruck05eer}. This has been extended to exponential integrators in \cite{arara19sbs} and more
recently in \cite{yang22anc}.

Consider a semi-linear non-autonomous SDE of the form
\begin{equation}
  \label{eq:semilinearSDEt}
  \dX(t) = A(t)X(t)\dt + \sum_{\diffusionindex=0}^{\diffusionnum}
  g_{\diffusionindex}(X(t),t) \star \dW_{\diffusionindex}(t).
\end{equation}

{
\begin{example}
  \label{ex:example_semilinear}
  The SDE presented in \cref{ex:example} can be written in this form:
  \[
    \dop\begin{pmatrix} R(t) \\ V(t) \end{pmatrix} =
    \underbrace{\begin{pmatrix} 0 & 1 \\ 0 & -\alpha(t) \end{pmatrix}}_{A(t)}
    \begin{pmatrix} R(t) \\ V(t) \end{pmatrix} +
    \underbrace{\begin{pmatrix} 0 \\ f_d(R(t),t) \end{pmatrix}}_{g_0} \dt +
    \underbrace{\begin{pmatrix} 0 \\ f_s(R(t),V(t),t) \end{pmatrix}}_{g_1} \dW_1(t).
 \]
\end{example}
}

The following class of exponential integrators is considered:
\begin{subequations}  \label{eq:semilinearRKt}
\begin{align}
  H_i &= Z_{i0}(A;t_n,h) Y_n + \sum_{\diffusionindex=0}^{\diffusionnum} \sum_{j=1}^{\stagenum}
        Z_{ij}^{(\diffusionindex)}(A;t_n,h)g_m(H_j,t_n+c_jh), &i&=1,\dots,\stagenum, \\
  Y_{n+1} &= z_{0}(A;t_n,h) Y_n + \sum_{\diffusionindex=0}^{\diffusionnum}\sum_{i=1}^{\stagenum}
            z_{i}^{(\diffusionindex)}(A;t_n,h)g_\diffusionindex(H_i,t_n+c_ih){.}
\end{align}
\end{subequations}
{The Magnus type
methods proposed by Yang et al.\ \cite{yang21aco} for noncommutative SDEs \cref{eq:multilin}
belong to this class of exponential integrators. For these methods, construction of the random
coefficients is discussed in \cite{komori23ffm}.}

Yang et al.\ \cite{yang22anc} derived B--series for the exact and numerical solution of this problem, given that
the matrix function $A(t)$ is commutative,
by expanding exponentials of the form $\expo{\int_{t_n}^{t_n+h}{A(s)\ds}}$. In the following, it will be demonstrated  how these series can be derived significantly simpler by the theory derived
above, and  without the  assumptions on commutativity of $A(t)$.

Let us start with discussing the B--series for the exact solution of \cref{eq:semilinearSDEt}.
This is already a non-autonomous, vertically split system, so the results of \cref{sec:non-autonomous} hold.
%To simplify the notation, the partition index $\partitionindex$ is omitted. Further,
Let  $\bullet_\diffusionindex$ represent
$g_{\diffusionindex}$, $\diffusionindex = 0,1,\dotsc,\diffusionnum$, $\bullet_A$ represent the term $A{(t)}X{(t)}$  {that is linear} {with respect to $X{(t)}$} and $\bullet_t$ {be the node corresponding to} $t$. Clearly,
the term $A(t)X$ can only be differentiated once with respect to $X$, thus
$F(\tau=[\tau_1,\dotsc,\tau_{\kappa}]_{A})=0$ whenever {there is} more than one $\tau_k \in
\cup_{\diffusionindex=0}^{\diffusionnum} T_\diffusionindex \cup T_A$, that is, a node $\bullet_A$ can have maximum one child which is not in $T_t$.
The following result is thus a consequence of \cref{cor:autonomous}. %\cref{def:trees}, \cref{def:elementarydifferentials} and \cref{thm:exact}.
\begin{corollary} \label{cor:semilinear}
  The solution $X(t)$ of \cref{eq:semilinearSDEt} can be written as a B-series where
  \begin{enumerate}
    \item The set of trees  $T=T_g \cup T_A$ is defined by
      \begin{align*}
        &\bullet_\diffusionindex \in T_g \text{ for } \diffusionindex = 0,\dots,
        \diffusionnum \text{ and } \bullet_A \in T_A, \\
	&\tau = [\tau_1,\dotsc,\tau_{\kappa}]_{\diffusionindex} \in T_g \text{ for all
        } \diffusionindex=0,\dotsc,\diffusionnum \text{ if }
         \tau_1, \tau_2, \dotsc, \tau_{\kappa} \in T\cup \{\bullet_t\}, \\
	&\tau = [\tau_1, \bullet_t^{\kappa-1}]_A \in T_A \text{ if }
        \tau_1 \in T\cup\{\bullet_t\}.
      \end{align*}
    \item The elementary differentials are given by
      \begin{align*}
	& F(\bullet_{\diffusionindex})(x_0,t_0) = g_m(x_0,t_0), \qquad
	F(\bullet_A)(x_0,t_0) = A(t_0)x_0, \qquad F(\bullet_t)(x_0,t_0) = 1, \\
	& F(\tau = [\tau_1,\tau_2,\dotsm,\tau_{\kappa}]_{\diffusionindex})(x_0,t_0)
	= \left(D_\tau g_{\diffusionindex}(x_0, t_0)\right)
	\big( F(\tau_1)(x_0,t_0),\dotsc,F(\tau_{\kappa})(x_0,t_0) \big), \\
	& F\left(\tau = [\tau_1, \bullet_t^{\kappa-1}]_A\right)(x_0, t_0)=
          \begin{cases}
            A^{(\kappa)}(t_0)\cdot x_0 & \text{if } \tau_1 = \bullet_t, \\
            A^{(\kappa-1)}(t_0)\cdot F(\tau_1)(x_0,t_0) & \text{otherwise}.
          \end{cases}
    \end{align*}
\end{enumerate}
\end{corollary}

We will now move our attention to the B--series of the exponential RK--method given in
\cref{eq:semilinearRKt}. In the case of a constant $A$, the coefficients are usually represented as
some kind of functions of $Ah$ or approximations of those, and the order conditions are
found by the series expansions of these functions. In the present case, with a time-dependent $A$,
we will assume that the coefficients $Z_{ij}^{(m)}$, { $i=1,\dots,s$, $j=0,\dots,s$,} and $z_i^{(m)}$, { $i=0,\dots,s$,} can be written in terms of a
B--series

\begin{align} \label{eq:coeffseries}
  Z(A;t_0,h) x_0 &= \sum_{\tau \in \bar{T}_A} \alpha(\tau)\cdot Z(\tau)(h) \cdot
                            F(\tau)(x_0,t_0) = B(Z,x_0,t_0;h)
\end{align}
where $\bar{T}_A \subset T$ is the subset of $T$ consisting of trees with no $\bullet_m$ nodes, $\diffusionindex=0,\dots,\diffusionnum$.
\begin{theorem} \label{thm:bseries_erk}
  Assume that the method coefficients $Z_{i0}(A;t_0,h)$, $Z_{ij}^{(m)}(A;t_0,h)$, $z_0(A;t_0,h)$ and $z_i^{{(m)}}(A;t_0,h)$ all can be written as B--series of the form \cref{eq:coeffseries}, with $Z_{i0}(\emptyset)=z_0(\emptyset)\equiv1$.% and $Z_{ij}^{(\diffusionindex)}(\emptyset)(0)=z_i(\emptyset)(0)=0$.
   Then the stage values $H_i$ and the numerical solution after one step $Y_i$ can both be written as B--series
  \[ H_i = B(\Phi_i,x_0,t_0;h),\; i=1,\dotsc,\stagenum,  \qquad Y_1 = B(\Phi,x_0,t_0;h), \]
  where trees and elementary differentials are defined in \cref{cor:semilinear} and with
  \begin{align*}
    \Phi_i(\emptyset)(h) & = \Phi(\emptyset)(h)=1, \\
    \Phi_i(\bullet_m)(h) &= \sum_{j=1}^\stagenum Z_{ij}^{(m)}(\emptyset)(h),  \qquad \Phi_i(\bullet_t)(h) = c_ih, \qquad i=1,\dotsc,\stagenum,  \\
    \Phi(\bullet_m)(h) &= \sum_{i=1}^{\stagenum} z_i^{{(m)}}(\emptyset)(h),  \\
    \Phi_i(\tau)(h) & =
    \begin{cases}
      Z_{i0}(\tau)(h) & \text{ if } \tau \in \bar{T}_A , \\
      \sum_{j=1}^\stagenum Z_{ij}^{(m)}(\vartheta)(h) \cdot \prod_{k=1}^{\kappa} \Phi_j(\tau_k)(h) & \text{otherwise},
    \end{cases} \\
    \Phi(\tau)(h) & =
    \begin{cases}
      z_{0}(\tau)(h) & \text{ if } \tau \in \bar{T}_A,  \\
      \sum_{i=1}^\stagenum z_{i}^{(m)}(\vartheta)(h) \cdot \prod_{k=1}^{\kappa} \Phi_i(\tau_k)(h) & \text{otherwise},
    \end{cases}
  \end{align*}
  where $(\vartheta,\{\delta\})$ is the only non-zero element of $SP(\tau)$ for which $\vartheta \in \bar{T}_A $ and $\delta=[\tau_1,\dots,\tau_\kappa]_m \in T_g$.
\end{theorem}
\begin{example}\label{ex:exptree}
  Let
  \[
    \tau = [[[\bullet_t^2]_A,\bullet_0]_1,\bullet_t]_A = \texa{,}
  \]
  {in which the stochastic node with color $\diffusionindex=1$ is represented by a white node with index 1, the deterministic ones with shape 0 or $t$ are represented by black nodes with the shape as index, and the deterministic nodes with shape $A$ are represented by red nodes with index $A$.}
  This is the same tree as presented in \cref{ex:thmexact}, where the {mapping of the} corresponding node indices {is given by}
  $110 \rightarrow 0$, $111 \rightarrow 1$, $120 \rightarrow A$ and $210 \rightarrow t$. {The
  elementary differential $F(\tau)$ can in this case be directly derived from \cref{cor:semilinear},
  $F(\tau) = {\dot{A}}D_{\tau_{{1}}}g_1 \big({\ddot{A}}\cdot,g_0 \big)$ where {$\tau_1=[[\bullet_t^2]_A,\bullet_0]_1$,} $\dot{A}=\dop
  A/\dt$ and $\ddot{A} = \dop^2 A/\dt^2$. The weight
  function for the exact solution $\varphi(\tau)$ has not changed from \cref{ex:thmexact},
      \[{\varphi(\tau) = \frac13\int_0^h  \int_0^{s_1}  s_2^4 \star \dW_1(s_2)  s_1 \ds_1,}\]
  and the
  numerical weight function can be derived from \cref{thm:bseries_erk} {by noting that according to \cref{ex:sp}, the only element $(\vartheta,\{\delta\})$ in $SP(\tau)$ with $\vartheta\in\bar{T}_A$ and $\delta\in T_g$ is given by $(\vartheta,\{\delta\})=\Bigg(\treera, \bigg\{\treena\bigg\}\Bigg)=\Bigg(\treerasemlexp, \bigg\{\treenasemlexp\bigg\}\Bigg)$}:}
  \[
    \Phi(\tau) = \sum_{i=1}^{\stagenum} z_{i}^{(1)}([\bullet_t]_A)\cdot \Phi_i([\bullet_t^2]_A) \cdot \Phi_i(\bullet_0)
    = \sum_{i,j=1}^{\stagenum} z_{i}^{(1)}([\bullet_t]_A)  \cdot Z_{i0}([\bullet_t^2]_A) \cdot Z_{ij}^{(0)}(\emptyset){.}
  \]
\end{example}

\begin{proof}[Proof of \cref{thm:bseries_erk}]
From \cref{lem:function_of_b-series} the following can be deduced if $\phi(\emptyset)\equiv 1$:
\[
  g_m(B(\phi,x_0,t_0;h),t_0+ch) = B(\psi^{(m)}_{\phi},x_0,t_0;h)
\]
where
\[
  \psi^{(m)}_{\phi}(\tau)(h) =
  \begin{cases}
    \prod_{k=1}^{\kappa} \phi(\tau_k)(h) & \text{if } \tau=[\tau_1,\dots,\tau_\kappa]_m, \\
    0 & \text{otherwise}
  \end{cases}
\]
with the definition $\phi(\bullet_t)=ch$. Further, by \eqref{eq:coeffseries}, $Z(A;t_0,h) = \partial_2  B(Z,x_0,t_0;h)$, so \cref{bsa2:lem:Bfprime} applies to $Z(A;t_0,h) g_m(B(\phi,x_0,t_0;h),t_0+ch)$.

Assuming that $H_i = B(\Phi_i,x_0,t_0;h)$ and $Y_1=B(\Phi,x_0,t_0;h)$, the method \eqref{eq:semilinearRKt} can be written in terms of B--series
\begin{align*}
  B(\Phi_i,x_0,t_0;h) &= B(Z_{i,0},x_0,t_0;h) + \sum_{\diffusionindex=0}^\diffusionnum \sum_{j=1}^\stagenum \partial_2 B(Z_{ij}^{(m)},x_0,t_0;h) \cdot B(\psi^{(m)}_{\Phi_j},x_0,t_0;h), \\
  B(\Phi,x_0,t_0;h) &= B(z_{0},x_0,t_0;h) + \sum_{\diffusionindex=0}^\diffusionnum \sum_{i=1}^\stagenum \partial_2 B(z_{i}^{(m)},x_0,t_0;h) \cdot B(\psi^{(m)}_{\Phi_i},x_0,t_0;h).
\end{align*}

Now, by comparing term by term and applying \cref{bsa2:lem:Bfprime}, the statement of the theorem follows by induction over the height of the trees, starting with the induction hypothesis $\Phi(\emptyset)=\Phi_i(\emptyset)\equiv 1$. In addition, since each $\bullet_A$-node can have maximum one $\bullet_m$-node as root of a child tree, for each $\tau\in T$ there is maximum one non-zero element $(\vartheta,\{\delta\})$ in $SP(\tau)$ with $\vartheta\in\bar{T}_A$ and $\delta\in T_g$.
\end{proof}

We will conclude the paper with an example from Yang et.al.\ \cite{yang22anc}, in which low order terms are found and mean square convergence of order 1 for Stratonovich SDEs is proved. We do not intend to reproduce their results, but we will only use this example as a demonstration on how our theory can be applied.

\begin{example}\label{ex:expmidpointrule}
  Given a semi-linear Stratonovich SDE with one-dimensional noise:
  \[ \dX = A(t)X(t)\dt + g_0(X(t),t)\dt + g_1(X(t),t)\circ \dW.  \]
  The following exponential midpoint rule was proposed by Yang et.al.\ \cite{yang22anc}:
  \begin{align*}
    H_1 &= \expo{\int_{t_n}^{t_n+h/2} A(s) \ds}x_0 + \frac{h}{2} g_0(H_1,t_n+\frac{h}{2}) + \frac{\Delta W_n}{2}g_1(H_1, t_n+\frac{h}{2}), \\
    Y_{n+1} &= \expo{\int_{t_n}^{t_n+h}A(s) \ds}x_0 + \expo{\int_{t_n+h/2}^{t_n+h}A(s) \ds}\big(hg_0(H_1,t_n+\frac{h}{2}) + \Delta W_ng_1(H_1, t_n+\frac{h}{2})\big){,}
  \end{align*}
{where $\Delta W_n=W(t_{n+1})-W(t_n)$.}
Thus $\stagenum=1$ and $c_1=1/2$. The method coefficients are
\begin{align*}
  Z_{10}(A;t_n,h) &= \expo{\int_{t_n}^{t_n+h/2} A(s) \ds}  \\ & \hspace*{-1cm}=  I + \frac{h}{2}A_n + \frac{h^2}{8}({\dot{A}}_n+A_n^2) + \frac{h^3}{48}\big({\ddot{A}}_n+\frac{3}{2}({\dot{A}}_nA_n+A_n{\dot{A}}_n)+A_n^3\big) + \dotsc,   \\
  Z_{11}^{(0)}(A;t_n,h) &= \frac{h}{2}, \qquad Z_{11}^{(1)} = \frac{\Delta W_n}{2},  \\
  z_0(A;t_n,h) &= \expo{\int_{t_n}^{t_n+h} A(s) \ds} \\ &=  I + hA_n + \frac{h^2}{2}({\dot{A}}_n+A_n^2) + \frac{h^3}{6}\big({\ddot{A}}_n+\frac{3}{2}({\dot{A}}_nA_n+A_n{\dot{A}}_n)+A_n^3\big) + \dotsc,  \\
  z_1^{(0)}(A;t_n,h) &= \expo{\int_{t_n+h/2}^{t_n+h} A(s) \ds}h = hI + \frac{h^2}{2}A_n + \frac{h^3}{8}(3{\dot{A}}_n+A_n^2) + \dotsc,  \\
  z_1^{(1)}(A;t_n,h) &= \expo{\int_{t_n+h/2}^{t_n+h} A(s) \ds}\Delta W_n{I} = \Delta W_n + \frac{h\Delta W_n}{2}A_n + \frac{h^2\Delta W_n}{8}(3{\dot{A}}_n+A_n^2) + \dotsc,
\end{align*}

where $A_n=A(t_n)$, ${\dot{A}_n=\dot{A}(t_n)}$, ${\ddot{A}_n=\ddot{A}(t_n)}$.

It is clear that all method coefficients can be written as B--series of the form \cref{eq:coeffseries} and fulfill the conditions required in \cref{thm:bseries_erk}.

This can now be used to find the weight function $\Phi(\tau)(h)$ in \cref{ex:exptree}. As $\stagenum=1$ {and} $z_1^{(1)}([\bullet_t]_A) = 3h^2\Delta W_n/8$, $Z_{11}^{(0)}(\emptyset) = h/2$,
$Z_{10}([\bullet_t^2]_A) = h^3/48$ we get $\Phi(\tau) = 3h^6\Delta W /768$. Clearly
$\Phi(\tau)\not=\varphi(\tau)$, so the order condition related to this particular tree is not
satisfied, {which could also not be expected, as the order of the tree is $6.5$}.

Note that although only  SDEs with a commutative linear term are considered in \cite{yang22anc}, the above shows that this result also holds for the noncommutative case, as the Lie-bracket ${\dot{A}}A-A{\dot{A}}$ is a term of order $\mathcal{O}(h^3)$ and will thus have no effect on the order for a first order method.

\end{example}

\appendix
\section{\texorpdfstring{\cref{ex:thmexact}}{Example 2} applied to the SDE in \texorpdfstring{\cref{ex:example}}{Example 1}} \label{sec:appA}
We will here provide the details on how the result presented in \cref{ex:thmexact} for the elementary differential looks like more concrete for the SDE introduced in \cref{ex:example}.
The partitioned SDE in \cref{ex:example} is of the form
\begin{align*}
  dX^{(1)} &= g_0^{(1,1)}(X^{(1)},X^{(2)})\dt + g_0^{(1,2)}(X^{(1)},X^{(2)}) \dt + g_1^{(1,1)}(X^{(1)},X^{(2)}) \dW_1(t) \\
  dX^{(2)} &= g_0^{(2,1)}(X^{(1)},X^{(2)})\dt
\end{align*}
with
\[
  X^{(1)} = \begin{pmatrix} R \\ V \end{pmatrix}, \;
  X^{(2)} = t.
\]

The elementary differential $F(\tau)$ is then step by step %and the weight function $\varphi$ are
developed, starting with
\begin{align*}
F({\hspace{-2ex}\treenc})  &= g_0^{(1,1)} = \begin{pmatrix} 0 \\ f_d(R,t) \end{pmatrix},
  \qquad
  &F({\hspace{-2ex}\treenba}) &= g_0^{(1,2)} = \begin{pmatrix} V \\ -\alpha(t)V \end{pmatrix},
  \\
  F({\hspace{-2ex}\treenbb}) & = g_1^{(1,1)} = \begin{pmatrix} 0 \\ f_s(R,V,t) \end{pmatrix},
  &F(\hspace{-2ex}\treenb) &= g_0^{(2,1)} = 1
  .
\end{align*}
From \cref{def:elementarydifferentials} %and \cref{thm:exact}
 we derive
\begin{align*}
F(\treene)  &= \begin{pmatrix} 0 \\ -\ddot{\alpha}V \end{pmatrix},
  & F(\hspace{-5ex}\treena) &= \begin{pmatrix} 0 \\ \frac{\partial^2 f_s}{\partial V^2}\cdot (-\ddot{\alpha}  V)\cdot f_d \end{pmatrix},
\end{align*}
and finally for the tree
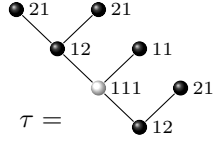
\[
\tau = \hspace{-8ex} \treeorig
\]
of \cref{ex:thmexact}
\[ F(\tau) = \begin{pmatrix} 0 \\ -\dot{\alpha}\cdot\frac{\partial^2 f_s}{\partial V^2}\cdot
                 (-\ddot{\alpha}V)\cdot f_d \end{pmatrix}.
\]
We notice that in the case when $f_s$ only depends on the position $R$, the elementary differential
for this tree is 0, thus even if the order condition is not satisfied for this particular tree, it
will still not contribute to the total error.
\section*{Acknowledgements}
The authors would like to thank the unknown referees for their very helpful comments.

\end{document}